\documentclass[11pt,draft]{article}
\usepackage{amssymb}
\usepackage{amsmath,amssymb}
\begin{document}
\title{{\bf Preimage entropy dimension of topological dynamical systems}}
\author{Lei Liu$^{1,2}$, Xiaomin Zhou$^{2}$, Xiaoyao Zhou$^{2*}$,  \\
\small 1 School of Mathematics and Information Science, Shangqiu Normal University,\\
    \small Shangqiu  476000, Henan P.R.China,\\
     \small  2 Department of Mathematics, University of Science and Technology of China,\\
      \small  Hefei, Anhui,230026, P.R.China\\
       \small    e-mail: mathliulei@163.com\\
        \small zxm12@mail.ustc.edu.cn\\
         \small    zhouxiaoyaodeyouxian@126.com,\\}
\date{}

\maketitle
\begin{quote}
\small{\bf Abstract.} We propose a new definition of preimage
entropy dimension for continuous maps on compact metric spaces,
investigate fundamental properties of the preimage entropy dimension,
and compare the preimage entropy dimension with the topological
entropy dimension. The defined preimage entropy dimension
holds various basic properties of topological entropy dimension,
for example, the preimage entropy dimension of a subsystem is
bounded by that of the original system and topologically conjugated
systems have the same preimage entropy dimension. Also, we
discuss the relation between the preimage entropy dimension
and the preimage entropy.

\end{quote}
\begin{quote}
\small{\bf MR Subject Classification (2010).} 54H20, 37B20.
\end{quote}
\begin{quote}
\end{quote}

\numberwithin{equation}{section}
\newtheorem{lm}{Lemma}[section]
\newtheorem{defi}{Definition}[section]
\newtheorem{prob}{Problem}
\newtheorem{thm}{Theorem}[section]
\newtheorem{pro}{Proposition}[section]
\newtheorem{exmp}{Example}[section]
\newtheorem{rmk}{Remark}
\newtheorem{cor}{Corollary}[section]
\newtheorem{con}{Conjecture}
\setlength\arraycolsep{2pt}
\baselineskip 16pt
\parskip 10pt


\section{Introduction}

\qquad In 1958, Kolmogorov applied the notion of entropy from
information theory to ergodic theory. Since then,
the concepts of entropies are useful for studying topological
and measure-theoretic structures of dynamical systems, that is,
topological entropy (see \cite{Adler-Konheim-McAndrew, Bowen1,
Bowen2}) and measure-theoretic entropy (see
\cite{Cornfeld-Fomin-Sinai-book, Kolmogorov-Tihomiorov}). For
instance, two conjugate systems have the same entropy and thus entropy
is a numerical invariant of the class of conjugated dynamical
systems. The theory of expansive dynamical systems has been closely
related to the theory of topological entropy \cite{Bowen-Walters,
Keynes-Sears, Thomas}. Entropy and chaos are closely related, for
example, a continuous map of interval is chaotic if and only if it
has a positive topological entropy \cite{Block-Coppel-book}.

In \cite{Hurley}, Hurley introduced several other entropy-like
invariants for noninvertible maps. One of these, which Nitecki and
Przytycki \cite{Nitecki-Przytycki} called pre-image branch entropy
(retaining Hurley's notation), distinguishing points according to
the branches of the inverse map. Cheng and Newhouse
\cite{Cheng-Newhouse} further extended the concept of topological
entropy of a continuous map and gave the concept of pre-image
entropy for compact dynamical systems. Several important pre-image
entropy invariant, such as pointwise pre-image, pointwise branch
entropy, partial pre-image entropy, and bundle-like pre-image
entropy, etc., have been introduced and their relationships with
topological entropy have been established.  The authors \cite{Cheng-Li2,MaKuangLi,MaWu,KuangChengLi} extended and studied as above some entropy-like
invariants for the non-autonomous discrete dynamical systems given
by a sequence of continuous self-maps of a compact topological space.

Although systems with positive entropy are much more complicated
than those with zero entropy, zero entropy systems have
various complexity, and have been studied many authors
(see \cite{Carvalho,Cheng-Li,Dou-Huang-Park1,Dou-Huang-Park2,
Ferenczi-Park,Huang-Yi,Huang-Park-Ye,Misiurewicz-Smital,Park}). These authors
adopted various methods to classify zero dynamical systems.
Carvalho \cite{Carvalho} introduced the notion of entropy dimension
to distinguish the zero topological entropy systems and obtained
some basic properties of entropy dimension. Cheng and Li further
discuss entropy dimension of the probabilistic and the topological
versions and gives a symbolic subspace to achieve zero topological
entropy, but with full entropy dimension.
Ferenczi and Park \cite{Ferenczi-Park} investigated a new entropy-like
invariant for the action of $\mathbb{Z}$ or $\mathbb{Z}^d$ on a
probability space.

In this paper we introduce the generalized $s-$preimage entropy and
preimage dimension of a topological dynamical systems, and study properties
of the preimage entropy dimension. These properties include that the
preimage entropy dimension of continuous maps on the finite spaces
is zero, the preimage entropy dimension of a contractive continuous map
is also zero and the preimage entropy dimension of a subsystem is bounded by that
of the original system. We prove that the $s-$preimage entropy of $T^n$ is
less than or equal to $n$ times the $s-$preimage entropy of $T$ and
the preimage entropy dimension of $T^n$ equals the preimage entropy
dimension of $T$. Main results show that $s-$preimage entropy is an invariant,
that is, topologically conjugate systems have the same $s-$preimage entropy and
the $s-$preimage entropy of a direct product is the sum of the $s-$preimage
entropies of the factors. Moreover, we discuss the relation between
preimag entropy dimension and preimage entropy, and obtain that
the preimage entropy dimension is less than or equal to one
if the preimage entropy is less than positive infinity,
the preimage entropy dimension is more than or equal to one
if the preimage entropy equals positive infinity and
the preimage entropy dimension equals one if
the preimage entropy is more than zero and it is less than positive infinity.


\section{Preimage entropy dimension}

\qquad A topological dynamical system  $(X,d,T)$ ($(X,T)$ for short) means that  $(X,d)$ is
a compact metric space together with  a continuous self-map $T:X\to X$.
Let $\mathbb{N}$ denote the set of all positive integers
and let $\mathbb{Z_+}=\mathbb{N}\cup$ $\{0\}$.

Let $n\in\mathbb{N}$. Define the metric $d_{T,n}$  on
$X$ by
\begin{eqnarray*}
d_{T,n}(x,y)=\max\limits_{0\leq j<n}d(T^j(x),T^j(y)).
\end{eqnarray*}
Given a
subset $K\subseteq X$,
a set $E\subseteq K$ is an $(n,\epsilon,K,T)$-separated set if, for
any $x\neq y$ in $E$, one has $d_{T,n}(x,y)>\epsilon$.  Let the quantity $r(n,\epsilon,K,T)$
 be the maximal cardinality of $(n,\epsilon,K,T)$-separated
subset of $K$. A subset $E\subseteq K$ is an
$(n,\epsilon,K,T)$-spanning set if, for every $x\in K$, there is a
$y\in E$ such that $d_{T,n}(x,y)\leq\epsilon$. Let
$s(n,\epsilon,K,T)$ be the minimal cardinality of any
$(n,\epsilon,K,T)$-spanning set.  It
is standard that for any subset $K\subseteq X$,
\begin{eqnarray*}
s(n,\epsilon,K,T)\leq r(n,\epsilon,K,T)\leq
s(n,\epsilon/2,K,T).
\end{eqnarray*}
\begin{defi}
{\rm\cite{Cheng-Newhouse}} Let $(X,T)$ be a  topological dynamical system.
Let $\epsilon>0$ and
$n\in\mathbb{N}$. Then
\begin{align*}
h_{pre}(T)=\lim\limits_{\epsilon\to
0}\limsup\limits_{n\to\infty}\frac{1}{n}\log\sup\limits_{x\in
X,k\geq n}r(n,\epsilon, T^{-k}(x),T)
\end{align*}
 is called the preimage
entropy of $T$.
 \label{defi2.1}
\end{defi}

Note that if $T$ is a homeomorphism, then $h_{pre}(T)=0$.
\begin{defi}
{\rm\cite{Cheng-Li}} Let $(X,T)$ be a topological dynamical system and
$s\in\mathbb{Z_+}$. Then $s$-topological entropy of $T$ is defined as
\begin{eqnarray*}
D(s,T)=\lim\limits_{\epsilon\to 0}\limsup\limits_{n\to\infty}\frac{1}{n^s}\log r(n,\epsilon,X,T)
=\lim\limits_{\epsilon\to 0}\limsup\limits_{n\to\infty}\frac{1}{n^s}\log s(n,\epsilon,X,T).
\end{eqnarray*}
 \label{defi2.2}
\end{defi}

In\cite{Carvalho}, the author proved that the $s$- topological entropy $D(s,T)$ shares the following property.

\begin{pro}\label{prop2.1}
\begin{description}
\item[(1)]
The map $s>0\mapsto D(s,T)$ is positive and decreasing with $s$;
\item[(2)]
There exists $s_0\in[0,+\infty]$, such that
\begin{eqnarray*}
D(s,T)=\left\{
\begin{array}{ll}
 +\infty,
&\mbox{\rm if} ~ 0<s<s_0,
\\
\;
 \\
 0,&\mbox{\rm if}~s>s_0;
~~~~~~~~~~~~~~~~~~~~~~~~~~~
\end{array}
\right.
\end{eqnarray*}
\item[(3)]
 $D(s,T^m)\leq m^sD(s,T)$.
\end{description}

\end{pro}
Proposition \ref{prop2.1}(2) indicates that the value of $D(s,T)$ jumps from infinity to $0$ at the two
sides of some point $s_0$, which is similar to a fractal measure. Analogous to the fractal
dimension, Cheng and Li \cite{Cheng-Li} defined the entropy dimension of $T$ as follows:
\begin{eqnarray*}
D(T)=\sup\{s>0:D(s,T)=\infty\}=\inf\{s>0:D(s,T)=0\}.
\end{eqnarray*}

Now, analogous to entropy dimension, we begin our process to introduce our new definition of
preimage entropy dimension.
\begin{defi}
Let $(X,d,T)$ be a topological dynamical system. Let $\epsilon>0$, $n\in\mathbb{N}$ and
$s>0$. Then $s$-preimage entropy of $T$ is defined as
\begin{eqnarray*}
D^d_{pre}(s,T)&=\lim\limits_{\epsilon\to
0}\limsup\limits_{n\to\infty}\frac{1}{n^s}\log\sup\limits_{x\in
X,k\geq n}r(n,\epsilon, T^{-k}(x),T)\\
&=\lim\limits_{\epsilon\to
0}\limsup\limits_{n\to\infty}\frac{1}{n^s}\log\sup\limits_{x\in
X,k\geq n}s(n,\epsilon, T^{-k}(x),T).
\end{eqnarray*}
The preimage  entropy dimension of $T$ is given by
\begin{eqnarray*}
D^d_{pre}(T)=\inf\{s>0:D^d_{pre}(s,T)=0\}=\sup\{s>0:D^d_{pre}(s,T)=\infty\}.
\end{eqnarray*}
 \label{defi2.3}
\end{defi}
When there is no confuse, we use $D_{pre}(s,T)$ instead of $D^d_{pre}(s,T)$ and $D_{pre}(T)$ instead of $D^d_{pre}(T).$
Clearly, $D_{pre}(1,T)=h_{pre}(T)$. If $T$ is a homeomorphism, then $D_{pre}(s,T)=0$, further,
$D_{pre}(T)=0$. When $X$ needs to be explicitly mentioned, we write $D_{pre}(T,X)$ instead of
$D_{pre}(T)$.

\begin{pro}
$D_{pre}(T)$ is independent of the choice of metric on $X$.
 \label{pro2.2}
\end{pro}
{\bf Proof.} We only prove that $D_{pre}(s,T)$ is independent
of the choice of metric on $X$. Let $d_1$
and $d_2$ be two compatible metrics on $X$. Then, by the compactness of $X$, for every $\epsilon>0$, there is a $\delta>0$ such
that, for all $x,y\in X$, if $d_1(x,y)<\delta$, then
$d_2(x,y)<\epsilon$. It follows that
$r(n,\epsilon,T^{-k}(x),T,d_2)\leq
r(n,\delta,T^{-k}(x),T,d_1)$ for all $x\in X$ and for every $n\in\mathbb{N}$ with $k\geq n$. This
shows that $D^{d_2}_{pre}(s,T)\leq
D^{d_1}_{pre}(s,T)$.  Interchanging $d_1$ and $d_2,$ this gives the
opposite inequality. Therefore,
$D^{d_1}_{pre}(s,T)=D^{d_2}_{pre}(s,T)$.  $\hfill{} \Box$
\begin{pro}
If $X$ is a finite set, then $D_{pre}(T)=0$.
 \label{pro2.3}
\end{pro}
{\bf Proof.} Let $\epsilon>0$, $n\in\mathbb{N}$ and $s>0$. Since $X$ is finite,
$r(n,\epsilon,T^{-k}(x),T)\leq card(X)$. Hence,
\begin{align*}
&D_{pre}(s,T)\\=&\lim\limits_{\epsilon\to
0}\limsup\limits_{n\to\infty}\frac{1}{n^s}\log\sup\limits_{x\in
X,k\geq n}r(n,\epsilon, T^{-k}(x),T)\\ \leq& \lim\limits_{\epsilon\to
0}\limsup\limits_{n\to\infty}\frac{1}{n^s}\log card(X)=0.
\end{align*}
 This shows
that $D_{pre}(s,T)=0$ for any $s>0$. Therefore, $D_{pre}(T)=0$.  $\hfill{} \Box$

We can easily obtain the following proposition by using the method of
\cite{Carvalho}.
\begin{pro}
\begin{description}
\item[(1)]
The map $s>0\mapsto D_{pre}(s,T)$ is positive and decreasing with $s$;
\item[(2)]
There exists $s_0\in[0,+\infty]$, such that
\begin{eqnarray*}
D_{pre}(s,T)=\left\{
\begin{array}{ll}
 +\infty,
&\mbox{\rm if} ~ 0<s<s_0,
\\
\;
 \\
 0,&\mbox{\rm if}~s>s_0.
~~~~~~~~~~~~~~~~~~~~~~~~~~~
\end{array}
\right.
\end{eqnarray*}
\end{description}
 \label{pro2.4}
\end{pro}

\begin{pro}
If $T:X\to X$ is a contractive continuous map, then $D_{pre}(T)=0$.
 \label{pro2.5}
\end{pro}
{\bf Proof.} If $T:X\to X$ is a contractive map, then $T$ diminishes
distance. Then
\begin{eqnarray*}
s(n,\epsilon, T^{-k}(x),T)\leq s(n-1,\epsilon, T^{-k}(x),T)\leq
\cdots\leq s(1,\epsilon, T^{-k}(x),T)
\end{eqnarray*}
and so
\begin{align*}
D_{pre}(s,T)=&\lim\limits_{\epsilon\to 0}
\limsup\limits_{n\to\infty}\frac{1}{n^s}\log\sup\limits_{x\in
X,k\geq n}s(n,\epsilon, T^{-k}(x),T)\\
\leq&\lim\limits_{\epsilon\to 0}
\limsup\limits_{n\to\infty}\frac{1}{n^s}\log\sup\limits_{x\in
X,k\geq n}s(1,\epsilon, T^{-k}(x),T)=0 ~{\rm for~all}~ s>0.
\end{align*}
This shows that $D_{pre}(s,T)=0$ for all $s>0$. Therefore,
$D_{pre}(T)=0$.  $\hfill{} \Box$

Denote by $K(X,T)$ the set of all $T$-invariant nonempty compact
subsets of $X$, that is, $K(X,T)=\{F\subseteq X:F\neq\emptyset, {\it
F}~{\rm is ~compact~and~}{\it T(F)\subseteq F}\}$. Since $X$ is
compact, it follows from $T(X)\subseteq X$ that
$K(X,T)\neq\emptyset$.

\begin{defi}
Let $(X,T)$ be a topological dynamical system. Let $\epsilon>0$ and $n\in\mathbb{N}$. For $F\in
K(X,T)$,
\begin{align*}
D_{pre}(s,T|_F,F)=\lim\limits_{\epsilon\to
0}\limsup\limits_{n\to\infty}\frac{1}{n^s}\log\sup\limits_{x\in
F,k\geq n}r(n,\epsilon,(T|_F)^{-k}(x),T|_F)
\end{align*}
is called
$s$-preimage entropy of $T$ on $F$, where $T|_F:F\to F$ is the induced
map of $T$, that is, for any $x\in F$, $T|_F(x)=T(x)$. The preimage entropy
dimension of $T$ restricted to $F$ is given by
\begin{eqnarray*}
D_{pre}(T|_F,F)=\inf\{s>0:D_{pre}(s,T|_F,F)=0\}=\sup\{s>0:D_{pre}(s,T|_F,F)=\infty\}.
\end{eqnarray*}
 \label{defi2.4}
\end{defi}

\begin{pro} \label{pro2.6}
If $F_1,F_2\in K(X,T)$ and $F_1\subseteq F_2$, then
$D_{pre}(T|_{F_1},F_1)\leq D_{pre}(T|_{F_2},F_2)$.
\end{pro}
{\bf Proof.} Let $\epsilon>0$ and $n,~k\in\mathbb{N}$ with $k\geq
n$, and let $x\in F_1$ and $E\subseteq (T|_{F_1})^{-k}(x)$ be an
$(n,\epsilon,(T|_{F_1})^{-k}(x),T|_{F_1})$-separated subset with
the maximal cardinality.   For $x\in F_1,$ we have
\begin{align*}
F_2\supseteq (T|_{F_2})^{-k}(x)\supseteq (T|_{F_1})^{-k}(x)\subseteq F_1\subseteq F_2.
\end{align*}
Hence, $E$ is an
$(n,\epsilon,(T|_{F_2})^{-k}(x),T|_{F_2})$-separated subset of
$(T|_{F_2})^{-k}(x)$. Therefore,
\begin{eqnarray*}
r(n,\epsilon,(T|_{F_1})^{-k}(x),T|_{F_1})\leq
r(n,\epsilon,(T|_{F_2})^{-k}(x),T|_{F_2}).
\end{eqnarray*}
Furthermore, we have
\begin{align*}
&\lim\limits_{\epsilon\to
0}\limsup\limits_{n\to\infty}\frac{1}{n^s}\log\sup\limits_{x\in
F_1,k\geq
n}r(n,\epsilon,(T|_{F_1})^{-k}(x),T|_{F_1})\\ \leq&
\lim\limits_{\epsilon\to
0}\limsup\limits_{n\to\infty}\frac{1}{n^s}\log\sup\limits_{x\in
F_2,k\geq n}r(n,\epsilon,(T|_{F_2})^{-k}(x),T|_{F_2}).
\end{align*}
This shows that $D_{pre}(s,T|_{F_1},F_1)\leq D_{pre}(s,T|_{F_2},F_2)$. By
Definition \ref{defi2.4}, we have
\begin{align*}
 D_{pre}(T|_{F_1},F_1)\leq D_{pre}(T|_{F_2},F_2).
\end{align*}  $\hfill{} \Box$

\begin{rmk}
Fix $s>0.$ For $F_i\in K(X,T),i\in\mathbb N,$ if $F=\bigcup\limits_{i=1}^{\infty}F_i$ and $F\in K(X,T)$,
then $D_{pre}(s,T|_F,F)\geq\sup\limits_{1\leq i<\infty}D_{pre}(s,T|_{F_i},F_i)$
and
$D_{pre}(T|_F,F)\leq D_{pre}(T,X)$.
\end{rmk}

\begin{exmp}
Let $(\Sigma_2,\sigma)$ be a one-sided symbolic
dynamical system, where $\Sigma_2=\{x=(x_n)_{n=0}^{\infty}: x_n\in
\{0,1\}~\mbox{for every}~ n\}$,
$\sigma(x_0,x_1,x_2,\cdots)=(x_1,x_2,\cdots)$. Then
$D_{pre}(\sigma)=1$.
\end{exmp}

Considering $\{0,1\}$ as a discrete space and putting product
topology on $\Sigma_2$, an admissible metric $\rho$ on
the space $\Sigma_2$ is defined by
\begin{eqnarray*}
\rho(x,y)=\sum\limits_{n=0}^{\infty}\frac{d(x_n,y_n)}{2^n},
\end{eqnarray*}
where \begin{eqnarray*}
d(x_n,y_n)=\left\{
\begin{array}{ll}
 0,
&\mbox{if} ~ x_n=y_n,
\\
\;
 \\
 1,&\mbox{if}~x_n\neq y_n,
~~~~~~~~~~~~~~~~~~~~~~~~~~~
\end{array}
\right.
\end{eqnarray*}
for $x=(x_0,x_1,\cdots),~y=(y_0,y_1,\cdots)\in\Sigma_2$.
By Robinson {\rm\cite{Robinson-book}},
$\Sigma_2$ is a compact metric space.

From Nitecki \cite{Nitecki} and Cheng-Newhouse \cite{Cheng-Newhouse},
$h_{pre}(\sigma)=\log 2$. Since
\begin{align*}
&\lim\limits_{\epsilon\to
0}\limsup\limits_{n\to\infty}\frac{1}{n^s}\log\sup\limits_{x\in
\Sigma_2,k\geq n}r(n,\epsilon, \sigma^{-k}(x),\sigma)\\
=&\lim\limits_{\epsilon\to
0}\limsup\limits_{n\to\infty}\frac{1}{n}\log\sup\limits_{x\in
\Sigma_2,k\geq n}s(n,\epsilon, \sigma^{-k}(x),\sigma)(n^{1-s})\\
=&\left\{
\begin{array}{ll}
 +\infty,
&\mbox{if} ~s<1,
\\
\;
 \log 2,&\mbox{if}~s=1,
 \\
\;
 0,&\mbox{if}~s>1.
~~~~~~~~~~~~~~~~~~~~~~~~~~~
\end{array}
\right.~~~~~~~~~~~~~~
\end{align*}
Therefore, $D_{pre}(\sigma)=1$.  $\hfill{} \Box$

\section{Main results}

\begin{thm}
Let $(X,T)$ be a topological dynamical system. Then
\begin{description}
\item[(1)]
if $s>0$ and $m\in\mathbb{N}$, then $D_{pre}(s,T^m)\leq m^sD_{pre}(s,T)$,
\item[(2)]
if $0<s\leq 1$ and $m\in\mathbb{N}$, then $D_{pre}(s,T^{m})\leq mD_{pre}(s,T)$,
\item[(3)]
if $m\in\mathbb{N}$, then $D_{pre}(T^m)=D_{pre}(T)$.
\end{description}
 \label{thm3.1}
\end{thm}
{\bf Proof.} {\bf(1)} Given $s>0$ and $m\in\mathbb{N}$.
Write $g=T^m$. Let $k\geq n$ and $x\in X$. It is clear that
\begin{eqnarray*}
r(n,\epsilon,g^{-k}(x),g)\leq r(mn,\epsilon,T^{-mk}(x),T).
\end{eqnarray*}
Hence, we have
\begin{align*}
&\limsup\limits_{n\to\infty}\frac{1}{n^s}\log\sup\limits_{x\in
X,k\geq n}r(n,\epsilon, g^{-k}(x),g) \\ \leq&\limsup\limits_{n\to\infty}\frac{1}{n^s}\log\sup\limits_{x\in
X,k\geq n}r(nm,\epsilon, T^{-mk}(x),T) \\
=&\limsup\limits_{n\to\infty}\frac{m^s}{(nm)^s}\log\sup\limits_{x\in X,k\geq n}r(nm,\epsilon, T^{-mk}(x),T)\\
\leq&\limsup\limits_{n\to\infty}\frac{m^s}{(nm)^s}\log\sup\limits_{x\in X,k\geq mn}r(nm,\epsilon, T^{-k}(x),T) \\
=&m^s\limsup\limits_{n\to\infty}\frac{1}{(nm)^s}\log\sup\limits_{x\in X,k\geq mn}r(nm,\epsilon, T^{-k}(x),T)\\
\leq& m^s\limsup\limits_{n\to\infty}\frac{1}{(nm)^s}\log\sup\limits_{x\in X,k\geq n}r(nm,\epsilon, T^{-k}(x),T).
\end{align*}
Furthermore,
\begin{align*}
 &\lim\limits_{\epsilon\to 0}\limsup\limits_{n\to\infty}\frac{1}{n^s}\log\sup\limits_{x\in
X,k\geq n}r(n,\epsilon, g^{-k}(x),g)\\ \leq &m^s\lim\limits_{\epsilon\to 0}\limsup\limits_{n\to\infty}\frac{1}{(nm)^s}\log\sup\limits_{x\in
X,k\geq n}r(nm,\epsilon, T^{-k}(x),T).
\end{align*}
 Therefore, $D_{pre}(s,T^m)\leq m^sD_{pre}(s,T)$.

{\bf(2)} Since $m^s\leq m$ for $0<s\leq 1$ and $m\in\mathbb{N}$, it
follows that $D_{pre}(s,T^{m})\leq mD_{pre}(s,T)$ by the above result
of (1).

{\bf(3)} Suppose $D_{pre}(s,T)=0$ for some $s>0$. From the above (1),
we have $D_{pre}(s,T^m)=0$ for every $m\in\mathbb{N}$. Hence,
$\{s>0:D_{pre}(s,T)=0\}\subseteq\{s>0:D_{pre}(s,T^m)=0\}$, further,
$\inf\{s>0:D_{pre}(s,T^m)=0\}\leq\inf\{s>0:D_{pre}(s,T)=0\}$. This
shows that $D_{pre}(T^m)\leq D_{pre}(T)$.
Next, we prove that $D_{pre}(T^m)\geq D_{pre}(T)$. Since
\begin{eqnarray*}
D_{pre}(s,T)=\lim\limits_{\epsilon\to 0}\limsup\limits_{n\to\infty}
\frac{1}{n^s}\log\sup\limits_{x\in X,k\geq n}s(n,\epsilon, T^{-k}(x),T),
\end{eqnarray*}
then given $\alpha>0$, there exists $x_0\in X$ such that
\begin{eqnarray}
D_{pre}(s,T)-\alpha< \lim\limits_{\epsilon\to 0}\limsup\limits_{n\to\infty}
\frac{1}{n^s}\log\sup\limits_{k\geq n}s(n,\epsilon, T^{-k}(x_0),T).
\end{eqnarray}
Fix $m\in\mathbb{N}$, we have
\begin{eqnarray*}
D_{pre}(s,T^m)=\lim\limits_{\epsilon\to 0}\limsup\limits_{n\to\infty}
\frac{1}{n^s}\log\sup\limits_{y\in X,k\geq n}s(n,\epsilon, T^{-mk}(y),T^m).
\end{eqnarray*}
Taking $y=T^{mk-k}(x_0)$, we get
\begin{eqnarray}
D_{pre}(s,T^m)\geq\lim\limits_{\epsilon\to 0}\limsup\limits_{n\to\infty}
\frac{1}{n^s}\log\sup\limits_{k\geq n}s(n,\epsilon, T^{-k}(x_0),T^m).
\end{eqnarray}
As $X$ is a compact space, $T,T^2,\cdots,T^m$ are uniformly continuous
on $X$, given $\epsilon>0$, there exists $0<\delta<\epsilon$ such that
\begin{eqnarray*}
d(x,y)<\delta\Longrightarrow\max\limits_{0\leq i<m}d(T^i(x),T^i(y))<\epsilon.
\end{eqnarray*}
Then we get $s(mn,\epsilon,T^{-k}(x_0),T)\leq s(n,\delta,T^{-k}(x_0),T^m)$.
By (3.1) and (3.2), we have
\begin{align*}
D_{pre}(s,T^m)\geq&\lim\limits_{\delta\to 0}\limsup\limits_{n\to\infty}
\frac{1}{n^s}\log\sup\limits_{k\geq n}  s(n,\delta, T^{-k}(x_0),T^m)\\
\geq&\lim\limits_{\epsilon\to 0}\limsup\limits_{n\to\infty}
m^s\frac{1}{(mn)^s}\log\sup\limits_{k\geq n}s(mn,\epsilon, T^{-k}(x_0),T)\\
\geq& m^s(D_{pre}(s,T)-\alpha).
\end{align*}
Letting $\alpha\to 0,$ this implies that $D_{pre}(s,T^m)\geq m^sD_{pre}(s,T)$ for all $s>0$.
This shows that $D_{pre}(T^m)\geq D_{pre}(T)$.  $\hfill{} \Box$

\begin{thm}
Let $(X,T)$ be a topological dynamical system. Then
\begin{description}
\item[(1)]
If $h_{pre}(T)<+\infty$, then $D_{pre}(T)\leq 1$.
\item[(2)]
If $h_{pre}(T)=+\infty$, then $D_{pre}(T)\geq 1$.
\item[(3)]
If $0<h_{pre}(T)<+\infty$, then $D_{pre}(T)=1$.
\end{description}
 \label{thm3.2}
\end{thm}
{\bf Proof.} {\bf(1)} If  $h_{pre}(T)<+\infty$, then
\begin{eqnarray*}
\lim\limits_{\epsilon\to
0}\limsup\limits_{n\to\infty}\frac{1}{n}\log\sup\limits_{x\in
X,k\geq n}r(n,\epsilon, T^{-k}(x),T)<+\infty.
\end{eqnarray*}
Therefore, for all $s>1$, we have
\begin{align*}
&\lim\limits_{\epsilon\to
0}\limsup\limits_{n\to\infty}\frac{1}{n^s}\log\sup\limits_{x\in
X,k\geq n}r(n,\epsilon, T^{-k}(x),T)\\=&\lim\limits_{\epsilon\to
0}\limsup\limits_{n\to\infty}\frac{1}{n^{s-1}}\cdot\frac{1}{n}\log\sup\limits_{x\in
X,k\geq n}r(n,\epsilon, T^{-k}(x),T)\\
=&\lim\limits_{n\to\infty}\frac{1}{n^{s-1}}\lim\limits_{\epsilon\to
0}\limsup\limits_{n\to\infty}\frac{1}{n}\log\sup\limits_{x\in
X,k\geq n}r(n,\epsilon, T^{-k}(x),T)\\
=&0\cdot h_{pre}(T)=0.
\end{align*}
Therefore, $D_{pre}(T)\leq 1$.

{\bf(2)} If $h_{pre}(T)=+\infty$, then for all $s<1$, we have
\begin{align*}
&\lim\limits_{\epsilon\to
0}\limsup\limits_{n\to\infty}\frac{1}{n^s}\log\sup\limits_{x\in
X,k\geq n}r(n,\epsilon, T^{-k}(x),T)\\=&\lim\limits_{\epsilon\to
0}\limsup\limits_{n\to\infty}\frac{1}{n^{s-1}}\cdot\frac{1}{n}\log\sup\limits_{x\in
X,k\geq n}r(n,\epsilon, T^{-k}(x),T)
=+\infty,
\end{align*}
which implies that $D_{pre}(T)\geq 1$.

{\bf(3)} If $h_{pre}(T)<\infty$, then we get from (1) that $D_{pre}(T)\leq 1$.
As $h_{pre}(T)>0$, it is at $s=1$  that the map $s\longmapsto D_{pre}(s,T)$
changes its value:
\begin{align*}
D_{pre}(s,T)=
\left\{
  \begin{array}{ll}
    +\infty, & s<1; \\
   h_{pre}(T), & s=1; \\
    0, & s>1.
  \end{array}
\right.
\end{align*}
This means that $D_{pre}(T)=1$.  $\hfill{} \Box$

Let $(X,T_1)$ and $(Y,T_2)$ be two topological dynamical systems. Then,
$(X,T_1)$ is an extension of $(Y,T_2)$, or $(Y,T_2)$ is a factor of
$(X,T_1)$ if there exists a surjective continuous map $\pi: X\to Y$
(called a factor map) such that $\pi\circ T_1(x)=T_2\circ\pi(x)$ for
every $x\in X$. If   $\pi$ is a homeomorphism, then $(X,T_1)$
and $(Y,T_2)$ are said to be topologically conjugate and the
homeomorphism $\pi$ is called a conjugate map.

\begin{thm}
Let $(X,d,T_1)$ and $(Y,d',T_2)$ be two topological dynamical systems.
If $(X,T_1)$ and $(Y,T_2)$ are topologically conjugate with a conjugate
map $\pi:X\to Y$, then $D_{pre}(s,T_1)=D_{pre}(s,T_2)$ for all $s>0$.
 \label{thm3.3}
\end{thm}
{\bf Proof.} Let $n\in\mathbb{N}$ and $k\geq n$.
Since $X$ is compact and $\pi$ is continuous, $\pi$ is
uniform continuous. So for any given $\epsilon>0$, there exists
$\delta>0$ such that $d(x_1,x_2)\geq\delta$ whenever
$d'(\pi(x_1),\pi(x_2))\geq\epsilon$. Fix $y\in Y$, let
$E(n,\epsilon,T_2^{-k}(y),T_2)\subseteq T_2^{-k}(y)$ be a maximal
$(n,\epsilon,T_2^{-k}(y),T_2)$-separated set for $T_2^{-k}(y)$, i.e.,
\begin{align*}
 card (E(n,\epsilon,T_2^{-k}(y),T_2))=r(n,\epsilon,T_2^{-k}(y),T_2).
\end{align*}
Let $G\subseteq X$ be a set by taking $x'=\pi^{-1}(y')$ for each
$y'\in E(n,\epsilon,T_2^{-k}(y),T_2)$. Then $card(G)=card (E(n,\epsilon,T_2^{-k}(y),T_2))$.
For the above $y$, set $x=\pi^{-1}(y)$. Since for all $x'\in G$, there exists only
$y'\in E(n,\epsilon,T_2^{-k}(y),T_2)$ such that $x'=\pi^{-1}(y')$. Hence,
$y=T_2^k(y')=T_2^k(\pi(x'))=\pi (T_1^k(x'))$, further, $T_1^k(x')=x$, that is,
$x'\in T_1^{-k}(x)$, which implies that $G\subseteq T_1^{-k}(x)$.
We claim that $G$ is a $(n,\delta,T_1^{-k}(x),T_1)$-separated set for $T_1^{-k}(x).$ In fact, for any $x_1,x_2\in G$, $\pi(x_1),\pi(x_2)\in
E(n,\epsilon,T_2^{-k}(y),T_2)$. Thus
\begin{eqnarray*}
d_n'(\pi(x_1),\pi(x_2)):=\max\limits_{0\leq i<n}\{d'(T_2^i\pi(x_1),T_2^i\pi(x_2))\}\geq\epsilon.
\end{eqnarray*}
that is, there exists $0\leq i_0<n$ such that $d'(T_2^{i_0}\pi(x_1),T_2^{i_0}\pi(x_2))\geq\epsilon$.
Since $\pi\circ T_1^{i_0}=T_2^{i_0}\circ\pi$, then $d(T_1^{i_0}(x_1),T_1^{i_0}(x_2))\geq\delta$
by the uniform continuity of $\pi$. Thus
\begin{eqnarray*}
d_n(x_1,x_2):=\max\limits_{0\leq i<n}\{d(T_1^i(x_1),T_1^i(x_2))\}\geq\delta.
\end{eqnarray*}
Therefore, $r(n,\delta,T_1^{-k}(x),T_1)\geq card(G)=card(E(n,\epsilon,T_2^{-k}(y),T_2))=r(n,\epsilon,T_2^{-k}(y),T_2)$,
that is, $r(n,\delta,T_1^{-k}(x),T_1)\geq r(n,\epsilon,T_2^{-k}(y),T_2)$. Furthermore, we have
\begin{align*}
&\lim\limits_{\delta\to 0}
\limsup\limits_{n\to\infty}\frac{1}{n^s}\log\sup\limits_{x\in
X,k\geq n}r(n,\delta, T_1^{-k}(x),T_1)
\\\geq&\lim\limits_{\epsilon\to 0}
\limsup\limits_{n\to\infty}\frac{1}{n^s}\log\sup\limits_{y\in
Y,k\geq n}r(n,\epsilon, T_2^{-k}(y),T_2).
\end{align*}
This shows that $D_{pre}(s,T_1)\geq D_{pre}(s,T_2)$ for all $s>0$.
Since $\pi:X\to Y$ is a conjugate map, then $\pi^{-1}:Y\to X$
is also a conjugate map for $(Y,T_2)$ and $(X,T_1)$. Similarly,  we have $D_{pre}(s,T_2)\geq D_{pre}(s,T_1)$. Therefore,
$D_{pre}(s,T_1)=D_{pre}(s,T_2)$.  $\hfill{} \Box$

Let $(X,d_1,T_1)$ and $(Y,d_2,T_2)$ be two topological dynamical systems. For the product space
$X\times Y$, define a map $T_1\times T_2:X\times Y\to X\times Y$ by
$(T_1\times T_2)(x,y)=(T_1(x),T_2(y))$. This map $T_1\times T_2$ is continuous
and $(X\times Y,T_1\times T_2)$ forms a topological dynamical system.
The metric $d$ on $X\times Y$ is given by
\begin{eqnarray*}
d((x_1,y_1),(x_2,y_2))=\max\{d_1(x_1,x_2),d_2(y_1,y_2)\}~{\rm for~any}~
(x_1,y_1),(x_2,y_2)\in X\times Y.
\end{eqnarray*}

\begin{thm}
Let $(X,d_1,T_1)$ and $(Y,d_2,T_2)$ be two topological dynamical systems. Then
\begin{align*}
D_{pre}(s,T_1\times T_2)=D_{pre}(s,T_1)+D_{pre}(s,T_2).
\end{align*}
 \label{thm3.4}
\end{thm}
{\bf Proof.} Let $\epsilon>0,n\in \mathbb{N}$ and $k\geq n$. Since balls in the product metric $d$
are products of balls on $X$ and $Y$, the same is true for balls in the
metric $d_n$. Hence, for $(x,y)\in X\times Y,$
\begin{eqnarray*}
s(n,\epsilon,(T_1\times T_2)^{-k}(x,y),T_1\times T_2)\leq
s(n,\epsilon,T_1^{-k}(x),T_1)\cdot s(n,\epsilon,T_2^{-k}(y),T_2),
\end{eqnarray*}
which implies that $D_{pre}(s,T_1\times T_2)\leq D_{pre}(s,T_1)+D_{pre}(s,T_2)$ for any $s>0$.
On the other hand, the product of any $(n,\epsilon,T_1^{-k}(x),T_1)$-separated set in $T_1^{-k}(x)$ for $T_1$
and any $(n,\epsilon,T_2^{-k}(y),T_2)$-separated set in $T_2^{-k}(y)$ for $T_2$ is an $(n,\epsilon,T_1^{-k}\times T_2^{-k}(x,y),T_1\times T_2)$-separated set
in $T_1^{-k}(x)\times T_2^{-k}(y)$ for $T_1\times T_2$. Hence,
\begin{eqnarray*}
r(n,\epsilon,(T_1\times T_2)^{-k}(x,y),T_1\times T_2)\geq
r(n,\epsilon,T_1^{-k}(x),T_1)\cdot r(n,\epsilon,T_2^{-k}(y),T_2),
\end{eqnarray*}
which implies that $D_{pre}(s,T_1\times T_2)\geq D_{pre}(s,T_1)+D_{pre}(s,T_2)$ for any $s>0$.$\hfill{} \Box$

\subsection*{Acknowledgment}
The author would like to thank Dr. Tu Siming for useful discussion.
The work was  supported by  the Fundamental Research Funds for the Central Universities (grant No. WK0010000035).

\end{document}